\documentclass[reqno,9pt]{amsart}

\usepackage{comment}
\usepackage{wrapfig}
\usepackage[hang,small,bf]{caption}
\usepackage[subrefformat=parens]{subcaption}
\usepackage{comment}
\usepackage{indentfirst}
\usepackage[top=20mm, bottom=20mm, left=30mm, right=30mm]{geometry}
\usepackage{amscd}
\usepackage{amsmath,amsthm,amssymb}
\usepackage{color}
\usepackage[dvipdfmx]{graphicx}
\usepackage{amsfonts}
\usepackage{cases}
\usepackage{bm}
\usepackage{pb-diagram}
\usepackage{braket}
\usepackage{multirow}
\usepackage{empheq}
\usepackage{tikz}

\usepackage{listings,jvlisting}

\usetikzlibrary{intersections,calc,arrows}

\numberwithin{equation}{subsection}
\theoremstyle{definition}
\newtheorem{defi}{Definition}[section]

\newtheorem{lem}[defi]{Lemma}

\newtheorem{theo}[defi]{Theorem}
\newtheorem{coro}[defi]{Corollary}

\newtheorem{conj}[defi]{Conjecture}

\newcommand{\Ker}{\mathrm{Ker}}

\newcommand{\im}{\mathrm{Im}}

\newcommand{\rank}{\mathrm{rank}}

\newcommand{\initial}{\mathrm{in}}

\newcommand{\pddim}{\mathrm{pd}}

\newcommand{\join}{\vee}
\newcommand{\meet}{\wedge}

\title{
Minimal free resolution of Crystal module
}

\author{Yohei Oshida}

\address{Graduate School of Engineering and Science, Shibaura Institute of Technology, 307 Minumaku Fukasaku, Saitama-City 337-8570, Japan.}

\email{mf20022@shibaura-it.ac.jp}


\begin{document}

\begin{abstract}
The join-meet ideal was introduced by Takayuki Hibi in 1987. It is binomial ideals that are defined by finite lattices. We study the join-meet ideal of non-distributive finite lattices that do not always satisfy modular. In particular, we work on the case of Crystal lattice which is one of them. It was introduced by Yohei Oshida in 2022. The Crystal module in the title is the residue class ring which is in relation to the join-meet ideal of Crystal lattice. In this paper, we give important inequalities about the minimal free resolution of the Crystal module. The important point about this result is that the projective dimension and Betti-number of the Crystal module can be evaluated by using the Koszul complex.
\end{abstract}

\maketitle{}

\setcounter{tocdepth}{2}
\tableofcontents

\section{Introduction}
Let $L$ be a finite lattice and $K[L]$ be the polynomial ring over a field $K$ whose variables are the elements of $L$. The ideal
\begin{eqnarray}
I_L=(\{ a b - (a \join b) (a \meet b) \mid a , b \in L \}) \subset K[L] \nonumber
\end{eqnarray}
is called the join-meet ideal of $L$. It was introduced in \cite{4} by Hibi. As shown by \cite{1} or \cite{4}, $L$ is distributive if and only if $I_{L}$ is a prime ideal. If $L$ is distributive lattices, $K[L]$ is called Hibi ring.

We study the join-meet ideal of non-distributive finite lattices that do not always satisfy modular. In particular, we work on the case of Crystal lattice, that is, the finite lattice $L_{k}(n_1,\cdots,n_ k)$ which is defined by
\begin{eqnarray}
L_{k}(n_1,\cdots,n_ k)=\{ s , x_{1,1} , \cdots , x_{1,n_1} , \cdots , x_{k,1} , \cdots , x_{k,n_k} , t \} \nonumber 
\end{eqnarray}
where $s < x_{i,1} < \cdots < x_{i,n_i} < t$ for $1 \leq i \leq k$. It was introduced  by Yohei Oshida \cite{2} in 2022.
We will not mention it at all here, but there is an important conjecture ''Crystal Conjecture'' with relation to the join-meet ideal of Crystal lattice; see \cite{2} in detail.

Now, by using inverse lexicographic order induced by
\begin{eqnarray}
s \prec x_{1,1} \prec x_{1,2} < \cdots \prec x_{1,n_1} \prec x_{2,1} \prec x_{2,2} \prec \cdots \prec x_{2,n_2} \prec  \cdots \prec x_{k,1} \prec x_{k,2} \prec \cdots < x_{k,n_k} \prec t, \nonumber
\end{eqnarray}
we have
\begin{eqnarray}
\initial_{\prec}(I_{L_{k}(n_1,\cdots,n_k)})=( \bigcup_{i=1}^k \bigcup_{j=1,j \neq i}^k \{ x_{i,r_1} x_{j,r_2} \mid 1 \leq r_1 \leq n_i , 1 \leq r_2 \leq n_j \} ). \nonumber
\end{eqnarray}
Note that $\prec$ is a compatible monomial order; see [1, Example 6.16] in detail. The residue class ring 
\begin{eqnarray}
K[L_{k}(n_1,\cdots,n_k)]/\initial_{<}(I_{L_{k}(n_1,\cdots,n_k)}) \nonumber
\end{eqnarray}
is called the \text{\bf{Crystal module}}. We recall grade Betti number by \cite{1} and \cite{3}. Let $R_{n_1,\cdots,n_k}=K[L_{k}(n_1,\cdots,n_k)]$. Let denote a minimal free resolution of $R_{n_1,\cdots,n_k}$-module $R_{n_1,\cdots,n_k}/\initial_{<}(I_{L_{k}(n_1,\cdots,n_k)})$ by
\begin{eqnarray}\label{main minimal}
&&
0
\leftarrow R_{n_1,\cdots,n_k}
\xleftarrow{\varphi^{(n_1,\cdots,n_k)}_{1}} 
F^{(n_1,\cdots,n_k)}_1
\xleftarrow{\varphi^{(n_1,\cdots,n_k)}_{2}} 
F^{(n_1,\cdots,n_k)}_2
\xleftarrow{\varphi^{(n_1,\cdots,n_k)}_{3}}
\cdots
\nonumber \\
&&
\cdots
\xleftarrow{\varphi^{(n_1,\cdots,n_k)}_{\ell_{n_1,\cdots,n_k}-1}}
F^{(n_1,\cdots,n_k)}_{\ell_{n_1,\cdots,n_k}-1}
\xleftarrow{\varphi^{(n_1,\cdots,n_k)}_{\ell_{n_1,\cdots,n_k}}}
F^{(n_1,\cdots,n_k)}_{\ell_{n_1,\cdots,n_k}}
\leftarrow
0,
\end{eqnarray}
where $F^{(n_1,\cdots,n_k)}_i$ is a finitely generated grade free $R_{n_1,\cdots,n_k}$-module for $1 \leq i \leq \ell_{n_1,\cdots,n_k}$. Then, for free module $F^{(n_1,\cdots,n_k)}_i$, we can write
\begin{eqnarray}
F^{(n_1,\cdots,n_k)}_i=\bigoplus_{j} R_{n_1,\cdots,n_k}(-X^{\bm{a}_{i,j}}), \nonumber
\end{eqnarray}
where $X=s x_{1,1} \cdots x_{1,n_1} x_{2,1} \cdots x_{2,n_2} \cdots x_{k,1} \cdots x_{k,n_k} t$ and $\bm{a}_{i,j}$ is the element of $\mathbb{Z}^{n_1+\cdots+n_k+2}$. Hence, since we denote $R_{n_1,\cdots,n_k}(-X^{\bm{a}_{i,j}})=R_{n_1,\cdots,n_k}(- |\bm{a}_{i,j}|)$ where $| \cdot |$ is $1$-norm, the minimal free resolution (\ref{main minimal}) can be rewritten as follows:
\begin{eqnarray}\label{minimal free resolution}
&&
0
\leftarrow R_{n_1,\cdots,n_k}
\leftarrow
\bigoplus_{j} R_{n_1,\cdots,n_k}(- |\bm{a}_{1,j}|)^{B_{1,j}(n_1,\cdots,n_k)}
\leftarrow
\bigoplus_{j} R_{n_1,\cdots,n_k}(- |\bm{a}_{2,j}|)^{B_{2,j}(n_1,\cdots,n_k)}
\leftarrow
\nonumber \\
&&
\cdots 
\leftarrow
\bigoplus_{j} R_{n_1,\cdots,n_k}(- |\bm{a}_{\ell_{n_1,\cdots,n_k}-1,j}|)^{B_{\ell_{n_1,\cdots,n_k}-1,j}(n_1,\cdots,n_k)}
\leftarrow
\bigoplus_{j} R_{n_1,\cdots,n_k}(- |\bm{a}_{\ell_{n_1,\cdots,n_k},j}|)^{B_{\ell_{n_1,\cdots,n_k},j}(n_1,\cdots,n_k)}
\leftarrow
0 \nonumber
\end{eqnarray}
We are mainly interested in computing the projective dimention
\begin{eqnarray}
\ell_{n_1,\cdots,n_k}=\pddim(R_{n_1,\cdots,n_k}/\initial_{<}(I_{L_{k}(n_1,\cdots,n_k)})) \nonumber
\end{eqnarray}
and $i$-th Betti number
\begin{eqnarray}
\sum_{j=0}^{\infty} B_{i,j}(n_1,\cdots,n_k). \nonumber
\end{eqnarray}
By \cite{6}, we have
\begin{eqnarray}
2 \leq \pddim(R_{n_1,n_2}/\initial_{<}(I_{L_{k}(n_1,n_2)})) \nonumber
\end{eqnarray}
and
\begin{eqnarray}
\sum_{i=0}^{\infty} B_{1,i}(n_1,1)=2n_1-1, \sum_{i=0}^{\infty} B_{2,i}(n_1,1)=n_1(n_1-1), \sum_{i=0}^{\infty} B_{1,i}(n_1,2)=3n_1, \sum_{i=0}^{\infty} B_{2,i}(n_1,2)=n_1(n_1+1)-1. \nonumber
\end{eqnarray}
The above result is very important for $k=2$. On the other hand, by using the above result, we cannot obtain a precise evaluation of the projective dimension and $i$-th Betti number of the Crystal module in the case of $k>2$. Moreover, in general, it is very difficult that us to obtain a precise evaluation.

In this paper, we give very important inequalities about the projective dimension and $i$-th Betti number of the Crystal module. Those results hold for any k. Moreover, the important point in them is that we can obtain them by using the Koszul complex.  In short, we can obtain very important inequalities by the following lemma.

\begin{lem}[The projective dimention of Koszul complex]\label{main} \it For positive integer $n_1$, the minimal free resolution of Koszul complex $K[x_{1,1},\cdots,x_{1,n_1}]/(x_{1,1},\cdots,x_{1,n_1})$ is
\begin{eqnarray}\label{Kreso}
0 \leftarrow S_{n_1} \leftarrow \bigoplus_{\bm{a} \in A(n_1,1)} S({X_{n_1}}^{\bm{a}}) \leftarrow \bigoplus_{\bm{a} \in A(n_1,2)} S({X_{n_1}}^{\bm{a}}) \leftarrow \cdots \leftarrow \bigoplus_{\bm{a} \in A(n_1,n_1)} S({X_{n_1}}^{\bm{a}}) \leftarrow 0,
\end{eqnarray}
where $X_{n_1}=x_{1,1} x_{1,2} \cdots x_{1,n_1}$, $S_{n_1}=K[x_{1,1},\cdots,x_{1,n_1}]$ and
\begin{eqnarray}
A(n_1,j)=\Biggr \{ ( a_1 , a_2 , \cdots , a_{n_1} ) \in \mathbb{Z}^{n_1}; \sum_{i=1}^{n_1} a_i = j , a_i \in \{ 0 , 1 \} \Biggr \}. \nonumber
\end{eqnarray}
for $j=1,2,\cdots,n_1$.
\end{lem}

By using the above lemma, we obtain the following main results in this paper.

\begin{theo}\label{theo1}\it For positive integer $k$ with $k>1$, we have
\begin{eqnarray}\label{theo1of1}
\max \{ n_1, n_2 \} \leq \pddim (R_{n_1,n_2}/\initial_{\prec}(I_{L_2(n_1,n_2)})) \leq n_1+n_2+1
\end{eqnarray}
and
\begin{eqnarray}\label{theo1of2}
{}_{n_1} C _i \leq \sum_{j=0}^{\infty} B_{i,j}(n_1,n_2)
\end{eqnarray}
for $i=1,2,\cdots,n_1$.
\end{theo}

By the above main result, we obtain the following corollary.

\begin{coro}\label{theo2}
\it For positive integer $k$ with $k>2$, we have
\begin{eqnarray}\label{theo2of1}
\max \{ n_1, \cdots, n_k\}
\leq 
\pddim (R_{n_1,\cdots,n_k}/\initial_{\prec}(I_{L_k(n_1,\cdots,n_k)})) 
\leq n_1+\cdots+n_k+2
\end{eqnarray}
and
\begin{eqnarray}\label{theo2of2}
{}_{n_1} C _i \leq \sum_{j=0}^{\infty} B_{i,j}(n_1,\cdots,n_k)
\end{eqnarray}
for $i=1,2,\cdots,n_1$.
\end{coro}

As we will see the proof of Theorem \ref{theo1} and Corollary \ref{theo2} later, the structure of the Crystal module can be observed by using the Koszul complex.

The structure of this paper is as follows. In Section 2, we prove Lemma \ref{main}. In Section 3, we prove Theorem \ref{theo1} by using Lemma \ref{main}. In Section 4, we prove Corollary \ref{theo2}. In Section 5, we introduce applications about the method of proof of  Theorem \ref{theo1}. In Section 6, we give a new conjecture about Theorem \ref{theo1}.

In the following, unless otherwise noted, the discussion of this paper will proceed with 
\begin{eqnarray}
\max \{ n_1, \cdots, n_k\}=n_1. \nonumber
\end{eqnarray}

\section{The proof of Lemma \ref{main}}

In this section, we prove Lemma \ref{main}.\\

At first, for $n_1=1$, we have the following minimal free resolution of $S_{1}/(x_{1,1})$:
\begin{eqnarray}
0 \leftarrow S_{1} \leftarrow S_{1}(x_{1,1})=\bigoplus_{\bm{a} \in A(1,1)} S_{1}({X_{1}}^{\bm{a}}) \leftarrow 0. \nonumber
\end{eqnarray}
The above minimal free resolution  equals (\ref{Kreso}) with $n_1=1$.\\

Second, for $n_1=1$, we have the following free resolution of $S_{2}/(x_{1,1},x_{1,2})$:
\begin{eqnarray}\label{Kreso1}
0 \leftarrow S_{2} \xleftarrow{g} \bigoplus_{i=1}^{2} S_{2}(x_i) \leftarrow \bigoplus_{\bm{a} \in A(2,2)} S_{2}({X_{2}}^{\bm{a}})=S_{2}(x_1 x_2) \leftarrow 0
\end{eqnarray}
Since we have
\begin{eqnarray}
\text{rank}(\bigoplus_{i=1}^{2} S_{2}(x_i))=\# \{ x_1 , x_2 \}=2, \nonumber
\end{eqnarray}we have $\Ker(g) \neq 0$. Thus, we obtain
\begin{eqnarray}
1 < \pddim(S_{2}/(x_{1,1},x_{1,2})). \nonumber
\end{eqnarray}
Hence, we have
\begin{eqnarray}
\pddim(S_{2}/(x_{1,1},x_{1,2}))=2. \nonumber
\end{eqnarray}
Therefore, the free resolution (\ref{Kreso1})  equals (\ref{Kreso}) with $n_1=2$.\\

Finally, for $n_1>1$, suppose that we have the minimal free resolution (\ref{Kreso}). Then, we have the following minimal free resoulutions:
\begin{eqnarray}
&& 0 \leftarrow S_{1} \leftarrow \bigoplus_{\bm{a} \in A(1,1)} S_{1}({X_{1}}^{\bm{a}}) \leftarrow 0, \nonumber \\
&& 0 \leftarrow S_{2} \leftarrow \bigoplus_{i=1}^{2} S_{2}(x_i) \leftarrow \bigoplus_{\bm{a} \in A(2,2)} S_{2}({X_{2}}^{\bm{a}}) \leftarrow 0, \nonumber \\
&& \quad \vdots \nonumber \\
&& 0 \leftarrow S_{n_1-1} \leftarrow \bigoplus_{i=1}^{n_1-1} S_{n_1}(x_i) \leftarrow \bigoplus_{\bm{a} \in A(n_1-1,2)} S_{n_1-1}({X_{n_1-1}}^{\bm{a}}) \leftarrow \cdots \leftarrow \bigoplus_{\bm{a} \in A(n_1-1,n_1-1)} S_{n_1-1}({X_{n_1-1}}^{\bm{a}}) \leftarrow 0. \nonumber
\end{eqnarray}
By using an inclusion map $f$ and $E_{i,j}=\bigoplus_{\bm{a} \in A(i,j)} S_{1}({X_{1}}^{\bm{a}})$, we obtain the following diagram:
\[
  \begin{CD}
     0 @<<< S_{1} @<<< E_{1,1} @<<< 0 \\
     @. @V{f}VV @V{f}VV \\
     0 @<<< S_{2} @<<< E_{2,1} @<<<  E_{2,2} @<<< 0 \\
     @. @V{f}VV @V{f}VV @V{f}VV \\
     0 @<<< S_{3} @<<< E_{3,1} @<<< E_{3,2} @<<<  E_{3,3} @<<< 0 \\
     @. @V{f}VV @V{f}VV @V{f}VV @V{f}VV \\
     @. \vdots @. \vdots @. \vdots @. \vdots @. @. \vdots \\
     @. @V{f}VV @V{f}VV @V{f}VV @V{f}VV @. @V{f}VV \\
     0 @<<< S_{n_1-1} @<<< E_{n_1-1,1} @<<< E_{n_1-1,2} @<<<  E_{n_1-1,3} @<<< \cdots @<<< E_{n_1-1,n_1-1} @<<< 0. \\
  \end{CD}
\]
Hence, we obtain the following diagram:
\[
  \begin{CD}
     0 @<<< S_{1} @<<< E_{1,1} @<<< 0 \\
     @. @V{f}VV @V{f}VV \\
     0 @<<< S_{2} @<<< E_{2,1} @<<<  E_{2,2} @<<< 0 \\
     @. @V{f}VV @V{f}VV @V{f}VV \\
     0 @<<< S_{3} @<<< E_{3,1} @<<< E_{3,2} @<<<  E_{3,3} @<<< 0 \\
     @. @V{f}VV @V{f}VV @V{f}VV @V{f}VV \\
     @. \vdots @. \vdots @. \vdots @. \vdots @. @. \vdots \\
     @. @V{f}VV @V{f}VV @V{f}VV @V{f}VV @. @V{f}VV \\
     0 @<<< S_{n_1-1} @<{g_{1,1}}<< E_{n_1-1,1} @<{g_{1,2}}<< E_{n_1-1,2} @<{g_{1,3}}<<  E_{n_1-1,3} @<{g_{1,4}}<< \cdots @<{g_{1,n_1-1}}<< E_{n_1-1,n_1-1} @<<< 0 \\
     @. @V{f}VV @V{f}VV @V{f}VV @V{f}VV @. @V{f}VV \\
     0 @<<< S_{n_1} @<{g_{2,1}}<< E_{n_1,1} @<{g_{2,2}}<< E_{n_1,2} @<{g_{2,3}}<<  E_{n_1,3} @<{g_{2,4}}<< \cdots @<{g_{2,n_1-1}}<< E_{n_1,n_1-1} @<{g_{2,n_1}}<< \cdots @<<< 0. \\
  \end{CD}
\]

By the above diagram, since we have $f(E_{n_1-1,i}) \subsetneq E_{n_1,i}$ for $i=1,2,\cdots,n_1-1$, we have
\begin{eqnarray}
f(\im(g_{1,i+1}))=f(\Ker(g_{1,i})) \subsetneq \Ker(g_{2,i})=\im(g_{2,i+1}) \nonumber
\end{eqnarray}
for $i=1,2,\cdots,n_1-2$.
Thus, it follows from this result we have
\begin{eqnarray}
n_1-1 \leq \pddim(S_{n_1}/(x_1,\cdots,x_{n_1})). \nonumber
\end{eqnarray}
Suppose that $n_1-1=(S_{n_1}/(x_1,\cdots,x_{n_1}))$. Then, we have
\begin{eqnarray}
0=f(0)=f(\Ker(g_{1,n_1-1})) \subsetneq \Ker(g_{2,n_1-1})=\im(0)=0. \nonumber
\end{eqnarray}
This is a contradiction. Thus, we have
\begin{eqnarray}
n_1-1 < \text{pd}(S_{n_1}/(x_1,\cdots,x_{n_1})). \nonumber
\end{eqnarray}
By [5, Theorem 4.6]. we obtain $\pddim(S_{n_1}/(x_1,\cdots,x_{n_1}))=n_1.$ Thus, the free resolution
\begin{eqnarray}\label{Kreso3}
0 \leftarrow S_{n_1} \leftarrow E_{n_1,1} \leftarrow E_{n_1,2} \leftarrow  E_{n_1,3} \leftarrow \cdots \leftarrow E_{n_1,n_1-1} \leftarrow E_{n_1,n_1} \leftarrow 0 \nonumber
\end{eqnarray}
is minimal. Hence, for positive integer $n_1$, the minimal free resolution of $K[x_{1,1},\cdots,x_{1,n_1}]/(x_{1,1},\cdots,x_{1,n_1})$ is (\ref{Kreso}).

\section{The proof of Theorem \ref{theo1}}

By Leema \ref{main} and using $E_i:=\bigoplus_{\bm{a} \in A(n_1,i)} S_{n_1}({X_{n_1}}^{\bm{a}})$, we have
\[
  \begin{CD}
     0 @<<< S_{n_1} @<{\varepsilon_0}<< E_{1} @<{\varepsilon_1}<< E_{2} @<{\varepsilon_2}<<
     \cdots
     @<{\varepsilon_{n_1-2}}<<  E_{n_1-1}  @<{\varepsilon_{n_1-1}}<< E_{n_1}  @<<< 0.
  \end{CD}
\]
Now, let's denote by $h$ the map from $E_i$ to $F^{(n_1,1)}_{i}$ by setting $h(E_i)=E_i x_{2,1}$. Then, we obtain the following result by using $h$.
\[
  \begin{CD}
     \Ker{(\varepsilon_{i-1})} @>{h}>> h(\Ker{(\varepsilon_{i-1})}) @>{\subset}>> \Ker{(\varphi^{(n_1,1)}_{i-1})} \\
  @|    @|  @| \\
     \im{(\varepsilon_i)}   @>{h}>> h(\im{(\varepsilon_i)}) @>{\subset}>> \im{(\varphi^{(n_1,1)}_i)} \\
  \end{CD}
\]
Thus, we obtain the following diagram which satisfies $\Ker{(\varepsilon_{i-1})}=\im{(\varepsilon_i)}$, $\Ker{(\varphi^{(n_1,1)}_{i-1})}=\im{(\varphi^{(n_1,1)}_i)}$ for $i=2,3,\cdots,n_1-1$:
\[
  \begin{CD}
     F^{(n_1,1)}_{i-1} @<{\varphi^{(n_1,1)}_{i-1}}<< F^{(n_1,1)}_{i} @<{\varphi^{(n_1,1)}_i}<< F^{(n_1,1)}_{i+1} \\
  @A{h}AA  @A{h}AA @A{h}AA \\
     E_{i-1} @<{\varepsilon_{i-1}}<< E_{i}   @<{\varepsilon_i}<<  E_{i+1}
  \end{CD}
\]
By the above diagram, we have the following diagram which satisfies $\Ker{(\varepsilon_{i-1})}=\im{(\varepsilon_i)}$, $\Ker{(\varphi^{(n_1,1)}_{i-1})}=\im{(\varphi^{(n_1,1)}_i)}$ for $i=2,3,\cdots,n_1-1$:
 \[
  \begin{CD}
     0 @<<< R_{n_1,1} @<{\varphi^{(n_1,1)}_0}<< F^{(n_1,1)}_{1} @<{\varphi^{(n_1,1)}_1}<< F^{(n_1,1)}_{2} @<{\varphi^{(n_1,1)}_2}<< \cdots
     @<{\varphi^{(n_1,1)}_{n_1-2}}<<  F^{(n_1,1)}_{n_1-1}  @<{\varphi^{(n_1,1)}_{n_1-1}}<< F^{(n_1,1)}_{n_1}  @<<<  \cdots @<<< 0 \\
    @. @AA{h}A  @AA{h}A  @AA{h}A @. @AA{h}A @AA{h}A \\
     0 @<<< S_{n_1} @<{\varepsilon_0}<< E_{1} @<{\varepsilon_1}<< E_{2} @<{\phi_2}<<
     \cdots
     @<{\varepsilon_{n_1-2}}<<  E_{n_1-1}  @<{\varepsilon_{n_1-1}}<< E_{n_1}  @<<< 0 
  \end{CD}
\]
Note that we have $h(E_1) \subsetneq F^{(n_1,1)}_{1}$; see [1] or [3] in detail.
Since we have the minimal free resolution (\ref{main minimal}),
we obtain the following result by using an inclusion map $f$:
\[
  \begin{CD}
     \Ker{(\varphi^{(n_1,1)}_{i-1})} @>{f}>> f(\Ker{(\varphi^{(n_1,1)}_{i-1})}) @>{\subset}>> \Ker{(\varphi^{(n_1,n_2)}_{i-1})} \\
  @|    @|  @| \\
     \im{(\varphi^{(n_1,1)}_i)}   @>{f}>>  f(\im{(\varphi^{(n_1,1)}_i)}) @>{\subset}>> \im{(\varphi^{(n_1,n_2)}_i)} \\
  \end{CD}
\]
Thus, we obtain the following diagram which satisfies $\Ker{(\varphi^{(n_1,1)}_i)}=\im{(\varphi^{(n_1,1)}_i)}$ and $\Ker{(\varphi^{(n_1,n_2)}_i)}=\im{(\varphi^{(n_1,n_2)}_i)}$ for $i=2,3,\cdots,n_1-1$:
\[
  \begin{CD}
     F^{(n_1,n_2)}_{i-1} @<{\varphi^{(n_1,n_2)}_{i-1}}<< F^{(n_1,n_2)}_{i} @<{\varphi^{(n_1,n_2)}_i}<< F^{(n_1,n_2)}_{i+1} \\
  @A{f}AA  @A{f}AA @A{f}AA \\
     F^{(n_1,1)}_{i-1} @<{\varphi^{(n_1,1)}_{i-1}}<< F^{(n_1,1)}_{i}   @<{\varphi^{(n_1,1)}_i}<<  F^{(n_1,1)}_{i+1}
  \end{CD}
\]

Hence, we obtain the following diagram by the above diagram:
 \[
  \begin{CD}
     0 @<<< R_{n_1,n_2} @<{\varphi^{(n_1,n_2)}_0}<< F^{(n_1,n_2)}_1 @<{\varphi^{(n_1,n_2)}_1}<< F^{(n_1,n_2)}_2 @<{\varphi^{(n_1,n_2)}_2}<< \cdots
     @<{\varphi^{(n_1,n_2)}_{n_1-2}}<<  F^{(n_1,n_2)}_{n_1-1}  @<{\varphi^{(n_1,n_2)}_{n_1-1}}<< F^{(n_1,n_2)}_{n_1}  @<<<
      \cdots \\
    @. @AA{f}A  @AA{f}A  @AA{f}A @. @AA{f}A @AA{f}A \\
     0 @<<< R_{n_1,1} @<{\varphi^{(n_1,1)}_0}<< F^{(n_1,1)}_{1} @<{\varphi^{(n_1,1)}_1}<< F^{(n_1,1)}_{2} @<{\varphi^{(n_1,1)}_2}<< \cdots
     @<{\varphi^{(n_1,1)}_{n_1-2}}<<  F^{(n_1,1)}_{n_1-1}  @<{\varphi^{(n_1,1)}_{n_1-1}}<< F^{(n_1,1)}_{n_1}  @<<< 
      \cdots \\
    @. @AA{h}A  @AA{h}A  @AA{h}A @. @AA{h}A @AA{h}A \\
     0 @<<< S_{n_1} @<{\varepsilon_0}<< E_{1} @<{\varepsilon_1}<< E_{2} @<{\varepsilon_2}<<
     \cdots
     @<{\varepsilon_{n_1-2}}<<  E_{n_1-1}  @<{\varepsilon_{n_1-1}}<< E_{n_1}  @<<< 0 
  \end{CD}
\]
By the above diagram, we obtain
\begin{eqnarray}\label{hana}
n_1 \leq \pddim(R_{n_1,n_2}/\initial_{\prec}(I_{L_2(n_1,n_2)})).
\end{eqnarray}
By using [2, the proof of Theorem 1.1], since we have
\begin{eqnarray}
\initial_{\prec}(I_{L_2(n_1,n_2)})=(\{ x_{1,i} x_{2,j} \mid 1 \leq i \leq n_1 , 1 \leq j \leq n_2 \} \cup \{ s x_{1,i} t \mid 2 \leq i \leq n_1 \} \cup \{ s x_{2,i} t \mid 2 \leq i \leq n_2 \}), \nonumber
\end{eqnarray}
it follows from [5, Theorem 3.1] that we obtain
\begin{eqnarray}\label{hana1}
\pddim(R_{n_1,n_2}/\initial_{\prec}(I_{L_2(n_1,n_2)})) \leq (n_1 + n_2 + 2) - 1= n_1 + n_2 + 1.
\end{eqnarray}
Hence, we obtain (\ref{theo1of1}) by  (\ref{hana}) and (\ref{hana1}). Moreover, for $i=1,2,\cdots,n_1$, since we have
\begin{eqnarray}\label{main number}
\rank(E_i)
=\rank \Biggr (\bigoplus_{\bm{a} \in A(n_1,i)} S({X_{n_1}}^{\bm{a}}) \Biggr )
=\# \{ X^{\bm{a}} ; \bm{a} \in A(n_1,i) \}
={}_{n_1} C _i,
\end{eqnarray}
we obtain (\ref{theo1of2}) by the following diagram:
\[
  \begin{CD}
     \Ker{(\varepsilon_{i-1})} @>{f \circ h}>> f \circ h(\Ker{(\varepsilon_{i-1})}) @>{\subset}>> \Ker{(\varphi^{(n_1,n_2)}_i)} \\
  @|    @|  @| \\
     \im{(\varepsilon_i)}   @>{f \circ h}>> f \circ h(\im{(\varepsilon_i)}) @>{\subset}>> \im{(\varphi^{(n_1,n_2)}_i)} \\
  \end{CD}
\]

\section{The proof of Corollary \ref{theo2}}

In this section, we prove Corollary \ref{theo2}. Let denote by $f$ an inclusion map. Let's denote by $h$ the map from $S_{n_1}$ to $R_{n_1,1}$ by setting $h(u)=u x_{2,1}$, where $u$ is the element of $S_{n_1}$. By the proof of Theorem \ref{theo2}, we have the following diagram which satisfies $\Ker{(\varepsilon_{i-1})}=\im{(\varepsilon_i)}$, $\Ker{(\varphi^{(n_1,1)}_{i-1})}=\im{(\varphi^{(n_1,1)}_i)}$ and $\Ker{(\varphi^{(n_1,n_2)}_{i-1})}=\im{(\varphi^{(n_1,n_2)}_i)}$ for $i=2,3,\cdots,n_1-1$:

 \[
  \begin{CD}
     0 @<<< R_{n_1,n_2} @<{\varphi^{(n_1,n_2)}_0}<< F^{(n_1,n_2)}_1 @<{\varphi^{(n_1,n_2)}_1}<< F^{(n_1,n_2)}_2 @<{\varphi^{(n_1,n_2)}_2}<< \cdots
     @<{\varphi^{(n_1,n_2)}_{n_1-2}}<<  F^{(n_1,n_2)}_{n_1-1}  @<{\varphi^{(n_1,n_2)}_{n_1-1}}<< F^{(n_1,n_2)}_{n_1}  @<<<
      \cdots \\
    @. @AA{f}A  @AA{f}A  @AA{f}A @. @AA{f}A @AA{f}A \\
     0 @<<< R_{n_1,1} @<{\varphi^{(n_1,1)}_0}<< F^{(n_1,1)}_{1} @<{\varphi^{(n_1,1)}_1}<< F^{(n_1,1)}_{2} @<{\varphi^{(n_1,1)}_2}<< \cdots
     @<{\varphi^{(n_1,1)}_{n_1-2}}<<  F^{(n_1,1)}_{n_1-1}  @<{\varphi^{(n_1,1)}_{n_1-1}}<< F^{(n_1,1)}_{n_1}  @<<< 
      \cdots \\
    @. @AA{h}A  @AA{h}A  @AA{h}A @. @AA{h}A @AA{h}A \\
     0 @<<< S_{n_1} @<{\varepsilon_0}<< E_{1} @<{\varepsilon_1}<< E_{2} @<{\varepsilon_2}<<
     \cdots
     @<{\varepsilon_{n_1-2}}<<  E_{n_1-1}  @<{\varepsilon_{n_1-1}}<< E_{n_1}  @<<< 0 
  \end{CD}
\]

Now, we have the following result.
\[
  \begin{CD}
     \Ker{(\varphi^{(n_1,n_2)}_{i-1})} @>{f}>> f(\Ker{(\varphi^{(n_1,n_2)}_{i-1})}) @>{\subset}>> \Ker{(\varphi^{(n_1,n_2,\cdots,n_k)}_{i-1})} \\
  @|    @|  @| \\
     \im{(\varphi^{(n_1,n_2)}_i)}   @>{f}>>  f(\im{(\varphi^{(n_1,n_2)}_i)}) @>{\subset}>> \im{(\varphi^{(n_1,n_2,\cdots,n_k)}_i)} \\
  \end{CD}
\]
Thus, we have the following diagram which satisfies $\Ker{(\varphi^{(n_1,n_2)}_{i-1})}=\im{(\varphi^{(n_1,n_2)}_i)}$, $\Ker{(\varphi^{(n_1,n_2,\cdots,n_k)}_{i-1})}=\im{(\varphi^{(n_1,n_2,\cdots,n_k)}_i)}$ for $i=1,2,\cdots,n_1-1$.
 \[
  \begin{CD}
     0 @<<< R_{n_1,\cdots,n_k} @<{\varphi^{(n_1,\cdots,n_k)}_0}<< F^{(n_1,\cdots,n_k)}_1 @<{\varphi^{(n_1,\cdots,n_k)}_1}<<  \cdots
     @<{\varphi^{(n_1,\cdots,n_k)}_{n_1-2}}<<  F^{(n_1,\cdots,n_k)}_{n_1-1}  @<{\varphi^{(n_1,\cdots,n_k)}_{n_1-1}}<< F^{(n_1,\cdots,n_k)}_{n_1}  @<<<
      \cdots \\
    @. @AA{f}A  @AA{f}A @. @AA{f}A @AA{f}A \\
     0 @<<< R_{n_1,n_2} @<{\varphi^{(n_1,n_2)}_0}<< F^{(n_1,n_2)}_1 @<{\varphi^{(n_1,n_2)}_1}<<  \cdots
     @<{\varphi^{(n_1,n_2)}_{n_1-2}}<<  F^{(n_1,n_2)}_{n_1-1}  @<{\varphi^{(n_1,n_2)}_{n_1-1}}<< F^{(n_1,n_2)}_{n_1}  @<<<
      \cdots
  \end{CD}
\]
Hence, we obtain the following diagram which satisfies $\Ker{(\varepsilon_{i-1})}=\im{(\varepsilon_i)}$, $\Ker{(\varphi^{(n_1,1)}_{i-1})}=\im{(\varphi^{(n_1,1)}_i)}$, $\Ker{(\varphi^{(n_1,n_2)}_{i-1})}=\im{(\varphi^{(n_1,n_2)}_i)}$ and $\Ker{(\varphi^{(n_1,n_2,\cdots,n_k)}_{i-1})}=\im{(\varphi^{(n_1,n_2,\cdots,n_k)}_i)}$ for $i=1,2,\cdots,n_1-1$ :
 \[
  \begin{CD}
     0 @<<< R_{n_1,\cdots,n_k} @<{\varphi^{(n_1,\cdots,n_k)}_0}<< F^{(n_1,\cdots,n_k)}_1 @<{\varphi^{(n_1,\cdots,n_k)}_1}<<  \cdots
     @<{\varphi^{(n_1,\cdots,n_k)}_{n_1-2}}<<  F^{(n_1,\cdots,n_k)}_{n_1-1}  @<{\varphi^{(n_1,\cdots,n_k)}_{n_1-1}}<< F^{(n_1,\cdots,n_k)}_{n_1}  @<<<
      \cdots \\
    @. @AA{f}A  @AA{f}A @. @AA{f}A @AA{f}A \\
     0 @<<< R_{n_1,n_2} @<{\varphi^{(n_1,n_2)}_0}<< F^{(n_1,n_2)}_1 @<{\varphi^{(n_1,n_2)}_1}<<  \cdots
     @<{\varphi^{(n_1,n_2)}_{n_1-2}}<<  F^{(n_1,n_2)}_{n_1-1}  @<{\varphi^{(n_1,n_2)}_{n_1-1}}<< F^{(n_1,n_2)}_{n_1}  @<<<
      \cdots \\
    @. @AA{f}A  @AA{f}A @. @AA{f}A @AA{f}A \\
     0 @<<< R_{n_1,1} @<{\varphi^{(n_1,1)}_0}<< F^{(n_1,1)}_{1} @<{\varphi^{(n_1,1)}_1}<< \cdots
     @<{\varphi^{(n_1,1)}_{n_1-2}}<<  F^{(n_1,1)}_{n_1-1}  @<{\varphi^{(n_1,1)}_{n_1-1}}<< F^{(n_1,1)}_{n_1}  @<<< 
      \cdots \\
    @. @AA{h}A  @AA{h}A @. @AA{h}A @AA{h}A \\
     0 @<<< S_{n_1} @<{\varepsilon_0}<< E_{1} @<{\varepsilon_1}<< 
     \cdots
     @<{\varepsilon_{n_1-2}}<<  E_{n_1-1}  @<{\varepsilon_{n_1-1}}<< E_{n_1}  @<<< 0 
  \end{CD}
\]

Therefore,  we obtain (\ref{theo2of1}) by the above diagram and [5, Theorem 4.6]. Moreover, we obtain (\ref{theo2of2}) by the above diagram and (\ref{main number}).

\section{Applications}

In this section, we introduce applications using the proof method of the Theorem \ref{theo1} and Corollary \ref{theo2}. This application is also deeply involved in Lemma \ref{main}.\\

Let denote $O_k(m_1,\cdots,m_k,M_1,\cdots,M_k)$ by the finite lattice 
\begin{eqnarray}
\{ s=t_0 , t_1 , t_2 , \cdots , t_{k-1} , t_k=t \} \cup \biggr ( \bigcup_{i=1}^{k} \{ z_{i,1} , z_{1,2} , \cdots , z_{i,m_i} \} \biggr ) \cup \biggr ( \bigcup_{i=1}^{k} \{ z_{i+k,1} , z_{i+k,2} , \cdots , z_{i+k,M_i} \biggr ), \nonumber
\end{eqnarray}
where $t_1 < t_2$, $t_2 < t_3$, $\cdots$, $t_{k-2}<t_{k-1}$ and  
\begin{eqnarray}
t_{i-1} < z_{i,1} < z_{i,2} < \cdots < z_{i,m_i} < t_i, \quad t_{i-1} < z_{i+k,1} < z_{i+k,2} < \cdots < z_{i+k,M_i} < t_i \nonumber 
\end{eqnarray}
for $i=1,2,\cdots,k$. By the above notation, we have
\begin{eqnarray}
\initial_{\prec}(I_{O_k(m_1,\cdots,m_k,M_1,\cdots,M_k)})=( \bigcup_{i=1}^k \{ x_{i,r_1} x_{i+k,r_2} \mid 1 \leq r_1 \leq m_i , 1 \leq r_2 \leq M_i \} ). \nonumber
\end{eqnarray}
where $\prec$ is a compatible monomial order ; see [1, 150 page] about definition of it. Then, we obtain the following theorem

\begin{theo}\label{theo4}For positive integer $k$ with $k>1$, we have
\begin{eqnarray}
\max \Biggr \{ m_1 , \cdots , m_k , M_1 , \cdots , M_k \Biggr \} \leq \pddim(K[O_k(m_1,\cdots,m_k,M_1,\cdots,M_k)]/\initial_{\prec}(I_{O_k(m_1,\cdots,m_k,M_1,\cdots,M_k)})) \leq k + 1 + \sum_{i=1}^k ( m_i + M_i ). \nonumber 
\end{eqnarray}
\end{theo}

In the following, unless otherwise noted, the discussion of this paper will proceed with 
\begin{eqnarray}
\max \Biggr \{ m_1 , \cdots , m_k , M_1 , \cdots , M_k \Biggr \}=m_1. \nonumber
\end{eqnarray}
Moreover, let denote the minimal free resolution of $K[O_k(m_1,\cdots,m_k,M_1,\cdots,M_k)]/\initial_{\prec}(I_{O_k(m_1,\cdots,m_k,M_1,\cdots,M_k)})$ by the following:
\begin{eqnarray}\label{main minimal}
&&
0
\leftarrow O_{m_1,\cdots,m_k,M_1,\cdots,M_k}
\xleftarrow{\phi^{(m_1,\cdots,m_k,M_1,\cdots,M_k)}_{1}} 
V^{(m_1,\cdots,m_k,M_1,\cdots,M_k)}_1
\xleftarrow{\phi^{(m_1,\cdots,m_k,M_1,\cdots,M_k)}_{2}} 
\cdots
\nonumber \\
&&
\cdots
\xleftarrow{\phi^{(m_1,\cdots,m_k,M_1,\cdots,M_k)}_{\ell_{m_1,\cdots,m_k,M_1,\cdots,M_k}-1}}
V^{(n_1,\cdots,n_k)}_{\ell_{n_1,\cdots,n_k}-1}
\xleftarrow{\phi^{(m_1,\cdots,m_k,M_1,\cdots,M_k)}_{\ell_{m_1,\cdots,m_k,M_1,\cdots,M_k}}}
V^{(m_1,\cdots,m_k,M_1,\cdots,M_k)}_{\ell_{m_1,\cdots,m_k,M_1,\cdots,M_k}}
\leftarrow
0,
\end{eqnarray}
where $V^{(m_1,\cdots,m_k,M_1,\cdots,M_k)}_i$ is a finitely generated grade free $O_{m_1,\cdots,m_k,M_1,\cdots,M_k}$-module for $1 \leq i \leq \ell_{m_1,\cdots,m_k,M_1,\cdots,M_k}$.\\




Let's prove Theorem \ref{theo4}.

\begin{proof}[The proof of Theorem \ref{theo4}]Let let denote $Z_i$ by $z_{1,1} z_{1,2} \cdots z_{1,m_1}$. By Lemma \ref{main}, for $i=1,2,\cdots,k$, we have the minimal free resolution
\begin{eqnarray}\label{ZZZ}
0 \leftarrow K[Z] \xleftarrow{\varepsilon_0} \bigoplus_{\bm{a} \in A(m_1,1)} K[Z](Z^{\bm{a}}) \xleftarrow{\varepsilon_1} \bigoplus_{\bm{a} \in A(m_1,2)} K[Z](Z^{\bm{a}}) \xleftarrow{\varepsilon_2} \cdots \xleftarrow{\varepsilon_{m-1}} \bigoplus_{\bm{a} \in A(m_1,m_1)} K[Z](Z^{\bm{a}}) \leftarrow 0,
\end{eqnarray}
where  $K[Z]:=K[z_{1,1}, z_{1,2}, \cdots, z_{1,m_1}]$ and
\begin{eqnarray}
A(m_1,j):=\Biggr \{ ( a_{1,1} , a_{1,2} , \cdots , a_{1,m_1}) \in \mathbb{Z}^{m_i}; \sum_{i=1}^{m_1} a_{1,i} = j , a_{1,i} \in \{ 0 , 1 \} \Biggr \}. \nonumber
\end{eqnarray}
for $j=1,2,\cdots,m_1$.\\

Let's denote $\bigoplus_{\bm{a} \in B(m_1,i)} K(Z^{\bm{a}})$ by $E_i$. Let's denote by $h$ the map from $E_i$ to $K[z_{1,1},z_{i,2},\cdots,z_{i,m_1},z_{1+k,M_1}]$ by setting $h(E_i)=E_i z_{1+k,M_1}$, where $u$ is the element of $K[Z]$. Then, we obtain the following diagram by using $h$.
\[
  \begin{CD}
     \Ker{(\varepsilon_{i-1})} @>{h}>> h(\Ker{(\varepsilon_{i-1})}) @>{\subset}>> \Ker{(\phi^{(m_1,0,\cdots,0=m_k,1)}_{i-1})} \\
  @|    @|  @| \\
     \im{(\varepsilon_i)}   @>{h}>> h(\im{(\varepsilon_i)}) @>{\subset}>> \im{(\phi^{(m_1,0,\cdots,0=m_k,1)}_i)} \\
  \end{CD}
\]
Moreover, we obtain the following diagram by using an inclusion map $f$:
\[
  \begin{CD}
     \Ker{(\phi^{(m_1,0,\cdots,0=m_k,1)}_{i-1})} @>{f}>> f(\Ker{(\phi^{(m_1,0,\cdots,0=m_k,1)}_{i-1})}) @>{\subset}>> \Ker{(\phi^{(m_1,\cdots,m_k,M_1,\cdots,M_k)}_{i-1})} \\
  @|    @|  @| \\
     \im{(\phi^{(m_1,0,\cdots,0=m_k,1)}_i)}   @>{f}>> f(\phi^{(m_1,0,\cdots,0=m_k,1)}_i) @>{\subset}>> \im{(\phi^{(m_1,\cdots,m_k,M_1,\cdots,M_k)}_i)} \\
  \end{CD}
\]
Thus, we obtain the following diagram: 
\[
  \begin{CD}
     \Ker{(\varepsilon_{i-1})} @>{f \circ h}>> f \circ h(\Ker{(\varepsilon_{i-1})}) @>{\subset}>> \Ker{(\phi^{(m_1,\cdots,m_k,M_1,\cdots,M_k)}_{i-1})} \\
  @|    @|  @| \\
     \im{(\varepsilon_i)}   @>{f \circ h}>> f \circ h(\im{(\varepsilon_i)}) @>{\subset}>> \im{(\phi^{(m_1,\cdots,m_k,M_1,\cdots,M_k)}_i)} \\
  \end{CD}
\]
By the above diagram, we obtain the following diagram which satisfies $\Ker{(\varepsilon_{i-1})}=\im{(\varepsilon_i)}$, $\Ker{(\phi^{(m_1,\cdots,m_k,M_1,\cdots,M_k)}_{i-1})}=\im{(\phi^{(m_1,\cdots,m_k,M_1,\cdots,M_k)}_i)}$ for $i=2,3,\cdots,m_1-1$:
\[
  \begin{CD}
     V^{(n_1,\cdots,n_k)}_{i-1} @<{\phi^{(m_1,\cdots,m_k,M_1,\cdots,M_k)}_{i-1}}<< V^{(n_1,\cdots,n_k)}_{i} @<{\phi^{(m_1,\cdots,m_k,M_1,\cdots,M_k)}_i}<< V^{(n_1,\cdots,n_k)}_{i+1} \\
  @A{f \circ h}AA  @A{f \circ h}AA @A{f \circ h}AA \\
     E_{i-1} @<{\varepsilon_{i-1}}<< E_{i}   @<{\varepsilon_i}<<  E_{i+1}
  \end{CD}
\]
Since we have the minimal free resolution (\ref{ZZZ}), it follows from the above diagram that we obtain
\begin{eqnarray}\label{hanahana}
m_1 \leq \pddim(K[O_k(m_1,\cdots,m_k,M_1,\cdots,M_k)]/\initial_{\prec}(I_{O_k(m_1,\cdots,m_k,M_1,\cdots,M_k)})).
\end{eqnarray}
Moreover, we obtain the following inequality by [5, Theorem 4.6]:
\begin{eqnarray}\label{hanahana1}
\pddim(K[O_k(m_1,\cdots,m_k,M_1,\cdots,M_k)]/\initial_{\prec}(I_{O_k(m_1,\cdots,m_k,M_1,\cdots,M_k)})) \leq k + 1 + \sum_{i=1}^k ( m_i + M_i ). 
\end{eqnarray}
Hence, we obtain the following inequality by (\ref{hanahana}) and (\ref{hanahana1}):
\begin{eqnarray}
m_1 \leq \pddim(K[O_k(m_1,\cdots,m_k,M_1,\cdots,M_k)]/\initial_{\prec}(I_{O_k(m_1,\cdots,m_k,M_1,\cdots,M_k)})) \leq k + 1 + \sum_{i=1}^k ( m_i + M_i ). \nonumber 
\end{eqnarray}

\end{proof}

\section{Conjecture}

At the end of this paper, we introduce the following conjecture.

\begin{conj}\it For positive integer $n_1$, we have $\pddim(R_{n_1,1}/\initial_{\prec}(I_{L_2(n_1,1)}))=n_1$. 
\end{conj}

\section*{Comment}

We comment a few. Please note that I am not a Shibaura Institute of Technology student after April 1st.\\

\Large

\end{document}